\documentclass[10pt]{amsart}

\newcommand{\la}{\lambda}
\newcommand{\al}{\alpha}
\newcommand{\be}{\beta}
\newcommand{\ga}{\gamma}

\newcommand{\ka}{K\"ahler}

\newcommand{\norm}[1]{\Vert #1\Vert}

\newcommand{\ov}{\overline}

\newcommand{\F}{\mathcal{F}}
\newcommand{\Ll}{\mathcal{L}}
 
\newcommand{\HA}{\mathcal{H}}
\newcommand{\JJ}{\mathcal{J}}
\newcommand{\KK}{\mathcal{K}}
\newcommand{\D}{\mathcal{D}} 
\newcommand{\T}{\mathcal{T}}     
\newcommand{\V}{\mathcal{V}}     
\newcommand{\f}{\varphi}
\newcommand{\e}{\varepsilon}

\newcommand{\CC}{\mathbb{C}}   
\newcommand{\HH}{\mathbb{H}}   
\newcommand{\RR}{\mathbb{R}}   
\newcommand{\ZZ}{\mathbb{Z}}

\newcommand{\HP}[1]{\mathbb{H}P^{#1}}

\numberwithin{equation}{section}

\newtheorem{te}{Theorem}[section]
\newtheorem{pr}{Proposition}[section]
\newtheorem{co}{Corollary}[section]
\newtheorem{lm}{Lemma}[section]

\theoremstyle{definition}
\newtheorem{de}{Definition}[section]    

\theoremstyle{remark}
\newtheorem{re}{Remark}[section]
\newtheorem{ex}{Example}[section]

\begin{document}

\title[Weyl structures on quaternionic manifolds]
{Weyl structures on quaternionic manifolds. A state of the art.}
\author{Liviu Ornea}
\date{}
\address{University of Bucharest, Faculty of Mathematics\newline
14 Academiei str., 
70109 Bucharest, Romania}
\email{lornea@imar.ro, lornea@gta.unibuc.ro}
\thanks{The author is a member of EDGE, Research Training Network
  HPRN-CT-2000-00101, supported by The European Human Potential Programme}
\maketitle

        This is a survey on quaternion Hermitian Weyl (locally 
conformally quaternion K\"ahler) and hyperhermitian Weyl (locally  
conformally hyperk\"ahler) manifolds. These geometries appear by 
requesting the compatibility of some quaternion Hermitian or hyperhermitian 
structure with a Weyl structure. The motivation for such a study is two-fold: 
it comes from the constantly growing interest in Weyl (and Einstein-Weyl)  
geometry and, on the other hand, from the necessity of understanding the 
existing classes of quaternion Hermitian manifolds.

        Various geometries are involved in the following discussion. The first 
sections give the minimal background on Weyl geometry, quaternion 
Hermitian geometry and $3$-Sasakian geometry. The reader is supposed familiar 
with Hermitian (K\"ahler and, if possible, locally conformally K\"ahler) and 
metric contact (mainly Sasakian) geometry.

        All manifolds and geometric objects on them are supposed differentiable 
of class $\mathcal{C}^\infty$. 
\section{Weyl structures}
We present here the necessary background concerning Weyl structures on conformal 
manifolds. We refer to \cite{F}, \cite{G1}, \cite{H} or to the most recent 
survey \cite{CP} for more details and physical interpretation (motivation) 
for Weyl and Einstein-Weyl geometry.

        Let $M$ be a $n$-dimensional, paracompact, smooth manifold, $n\geq 2$. 
A $\mathrm{CO}(n)\simeq \mathrm{O}(n)\times \RR_+$ structure on $M$ is equivalent with the 
giving of a conformal class $c$ of Riemannian metrics. The pair $(M,c)$ is 
a \emph{conformal manifold}. 

        For each metric $g\in c$ one can consider the Levi-Civita connection 
$\nabla^g$, but this will not be compatible with the conformal class. Instead, 
we shall work with  $\mathrm{CO}(n)$-connections. Precisely: 
\begin{de}
A \emph{Weyl connection} $D$ on a conformal manifold $(M,c)$ is a torsion-free 
connection which preserves the conformal class $c$. We say that $D$   
defines a Weyl structure on $(M,c)$ 
and $(M,c,D)$ is a Weyl manifold.
\end{de}
Preserving the conformal class means that for any $g\in c$, there exists a 
$1$-form $\theta_g$ (called the Higgs field) such that 
$$Dg=\theta_g\otimes g.$$ This formula is conformally invariant in  
the following sense: 
\begin{equation}\label{th}
\text{if}\; h=e^{f}g,\;\;  f\in \mathcal{C}^\infty(M),\; \text{then}\;\;  
\theta_h=\theta_g-df.
\end{equation}
        Conversely, if one starts with a fixed Riemannian metric $g$ on $M$ 
and a fixed $1$-form $\theta$ (with $T=\theta^\sharp$), the connection
$$D=\nabla^g-\frac{1}{2}\{\theta\otimes Id+Id\otimes\theta-g\otimes T\}$$
is a Weyl connection, preserving the conformal class of $g$. Clearly, 
$(g,\theta)$ and $(e^{f}, \theta-df)$ define the same Weyl structure. 

        On a Weyl manifold $(M,c,D)$, Weyl introduced the 
\emph{distance curvature function}, a $2$-form defined by $\Theta=d\theta_g$.  
By \eqref{th}, the definition does not depend on $g\in c$. If 
$\Theta=0$, the cohomology class $[\theta_g]\in H^1_{dR}(M)$ is independent 
on $g\in c$. A Weyl structure with $\Theta=0$ is called \emph{closed}. 

        All these geometric objects can be interpreted as sections in tensor bundles of 
the bundle of scalars of weight $1$, associated to the bundle of linear frames of $M$ 
\emph{via} the representation $GL(n,\RR)\ni A\mapsto \mid \det A\mid ^{1/n}$. 
\emph{E.g.} $c$ is a section of $S^2T^*M\otimes L^2$, $\theta $ is a connection form in 
$L$ whose curvature form is exactly the distance curvature function etc. 
This also motivates the terminology. We refer to \cite{G1} for a systematic treatment 
of this viewpoint.

        A fundamental result on Weyl structures is the following 
"co-close\-deness lemma":
\begin{te}\cite{G2}
Let $(M,c)$ be a compact, oriented, conformal manifold of dimension $>2$. 
For any Weyl structure $D$ preserving $c$, there exists a unique (up 
to homothety) $g_0\in c$ such that the associated Higgs field $\theta_{g_0}$ is 
$g_0$-coclosed.
\end{te}
The metric $g_0$ provided by the theorem is called the \emph{Gauduchon metric} 
of the Weyl structure. 

        In Weyl geometry, the good notion of Einstein manifold makes use 
of the Ricci tensor associated to the Weyl connection:
$$Ric^D=\frac{1}{2}\sum_{i=1}^n\{g(R^D(X,e_i)Y,e_i)-g(R^D(X,e_i)e_i,Y)\}$$
where $g\in c$ and $\{e_i\}$ is a local $g$-orthonormal frame. The scalar 
curvature of $D$ is then defined as the conformal $\mathrm{trace}$ of $Ric^D$. For 
each choice of a $g\in c$, $Scal^D$ is represented by $Scal^D_g=\mathrm{trace}_gRic^D$. 
The relations between $Ric^D$ and $Ric^{\nabla^g}$ and, correspondingly, 
between their scalar curvatures are:
\begin{gather}
\label{28}Ric^D=Ric^{\nabla^g}+\delta^g\theta\cdot g-(n-2)\{\nabla^g\theta+
\norm{\theta}^2_g\cdot g-\theta\otimes\theta\}.\\
\label{29}Scal^D_g=Scal^{\nabla^g}+2(n-1)\delta^g\theta-(n-1)(n-2)\norm{\theta}^2_g.
\end{gather}

\begin{de}
A Weyl structure is \emph{Einstein-Weyl} if the symmetric part of the 
Ricci tensor $Ric^D$ of the Weyl connection is proportional to one (hence 
any) metric of $c$.
\end{de}
For an Einstein-Weyl structure, one has
\begin{equation}\label{24}
Ric^D=\frac{1}{n}Scal^D_g\cdot g-\frac{n-2}{2}d\theta
\end{equation}
for any $g\in c$. 

        Note that for an Einstein-Weyl structure, the scalar curvature $Scal^D $ need not 
be constant (this means $D$-parallel w.r.t. to $D$ as a section of $L^{-2}$). 
But, if the Weyl connection is precisely the Levi-Civita connection of a 
metric in $c$ (in this case the Weyl structure is called \emph{exact}), then 
$Scal^D$ is constant. 

        Observe that for any Einstein-Weyl structure and any $g\in c$ one 
has the formula
\begin{equation}\label{27}
\begin{split}
\frac{1}{n}dScal^D_g&=\frac{2n}{n-2}d\delta^g\theta-
2(\delta^g\theta)\theta-2\delta^g\nabla^g\theta+\delta^gd\theta-\\
&-2\nabla^g_T\theta-(n-3)d\norm{\theta}^2_g.
\end{split}
\end{equation}
This follows from \eqref{24}, \eqref{28} and \eqref{29}. 

        If $g$ is the Gauduchon metric, $\delta^g\theta=0$ and \eqref{27} 
reduces to 
$$\frac{1}{n}dScal^D_g+
2\delta^g\nabla^g\theta-\delta^gd\theta+
+2\nabla^g_T\theta+(n-3)d\norm{\theta}^2_g=0.$$
Contracting here with $\theta$ yields
\begin{equation}\label{37}
D\theta=\frac{1}{2}d\theta
\end{equation}
This, together with the relation between $D$ and $\nabla^g$ prove the 
first statement of the following extremely important result (the second 
statement will be proved in a more particular situation):
\begin{te}\label{gg}\cite{G1}
Let $D$ be an Einstein-Weyl structure on a compact, oriented manifold 
$(M,c)$ of dimension $>2$. Let $g$ be the Gauduchon metric in $c$ 
associated to $D$ and $\theta$ the corresponding Higgs field. If the Weyl structure  
is closed, but not exact, then 

        1) $\theta$ is $\nabla^g$ parallel: $\nabla^g\theta=0$ (in particular, also 
$g$-harmonic).

        2) $Ric^D=0$.
\end{te}
Odd dimensional spheres and products of spheres $S^1\times S^{2n+1}$ admit 
Einstein-Weyl  structures (note that $S^1\times S^2$ and $S^1\times S^3$ 
can bear no Einstein metric, cf. \cite{Hi}). Further examples, with  
$Ric^D=0$, will be 
the compact quaternion Hermitian Weyl and hyperhermitian Weyl manifolds.

\section{Quaternionic Hermitian manifolds}
This section is devoted to the introduction of quaternion Hermitian geometry. 
The standard references are \cite{Sa},  \cite{Be}, \cite{Sw}, \cite{AM2}.
\begin{de}
Let $(M,g)$ be a $4n$-dimensional Riemannian manifold. Suppose $End (TM)$ 
has a rank $3$ subbundle $H$ with transition functions in $\mathrm{SO}(3)$, 
locally generated by orthogonal  almost complex 
structures 
$I_\al$, $\al=1,2,3$ satisfying the quaternionic relations. Precisely:
\begin{equation}\label{unu}
I_\al^2=-Id, \; I_\al I_\be=\varepsilon_{\al\be\ga}I_\ga, \; g(I_\al\cdot,
I_\al\cdot)=
g(\cdot,\cdot)\; \al,\be,\ga=1,2,3
\end{equation}
where $\varepsilon_{\al\be\ga}$ is $1$ (resp. $-1$) when $(\al\be\ga)$ is an 
even (resp. odd)  
permutation of $(123)$ (such a basis of $H$ is called admissible.) 
The triple $(M,g,H)$ is called a \emph{quaternionic 
Hermitian manifold} whose \emph{quaternionic  bundle} is $H$.
\end{de}
Any (local) or global section of $H$ is called \emph{compatible}, but in 
general, $H$ has no global section. A striking example is $\HP{n}$, 
the quaternionic projective space. The three canonical almost complex structures 
of $\HH^{n+1}$ induced by multiplication with the imaginary quaternionic 
units descend to only local almost complex structures on $\HP{n}$ generating 
the bundle $H$. The metric is the one projected by the flat one on $\HH^{n+1}$, 
\emph{i.e.} the Fubini-Study metric written in quaternionic coordinates. 
Note that $\HP{1}$ is diffeomorphic with $S^4$, hence cannot bear any almost complex 
structure. Consequently, no greater dimensional $\HP{n}$ can have an almost complex 
structure neither, because this would be induced on any quaternionic projective line 
$\HP{1}$, contradiction.

        This shows that the case when $H$ is trivial is of a special importance 
and motivates
\begin{de}
A quaternionic Hermitian manifold with trivial quaternionic bundle is called 
a \emph{hyperhermitian manifold}.
\end{de}
        In this terminology, an admissible basis of a quaternion Hermitian 
manifold 
is a local almost hyperhermitian structure. 

        For a hyperhermitian manifold we shall always fix a (global) 
basis of $H$ satisfying the quaternionic relations, so we shall regard it as 
a manifold endowed with three Hermitian structures $(g,I_\al )$ related by 
the identities \eqref{unu}
The simplest example is $\HH^{n}$. But we shall encounter other many examples. 

        The analogy with Hermitian geometry suggests imposing conditions 
of \ka\ type. Let $\nabla^g$ be the Levi-Civita connection of the metric $g$. 
\begin{de}
A quaternionic Hermitian manifold $(M,g,H)$ of dimension at least $8$ 
is \emph{quaternion \ka\ } if 
$\nabla^g$ parallelizes $H$, \emph{i.e.} $\nabla^gI_\al=a_\al^\be \otimes I_\be$ (with 
sqew-symmetric matrix of one forms $(a_\al^\be)$). 

        A hyperhermitian manifold is called \emph{hyperk\"ahler} if 
$\nabla^gI_\al=0$ for $\al=1,2,3.$
\end{de}
\begin{re}
This definition of quaternion K\"ahler manifold is redundant in dimension $4$. 
As S. Marchiafava proved (see \cite{Ma2}) that any four-dimensional isometric 
submanifold of a quaternion \ka\ 
manifold whose tangent bundle is invariant to each element of $H$ is Einstein and 
self-dual, one takes this as a definition. We won't be concerned with 
dimension $4$ in this report.
\end{re}
Note that, unless in the complex case, here the parallelism of $H$ does not 
imply the integrability of the single almost complex structures.

\begin{ex}
$\HH^{n}$ with its flat metric is hyperk\"ahler. By a result of A. Beauville, 
the $K^3$ surfaces also, 
see \cite{Be}, Chapter 14. The irreducible, symmetric  quaternion \ka\ were classified 
by J. Wolf. Apart $\HP{n}$, the compact ones are: the Grassmannian of oriented  
$4$-planes in $\RR^{m}$, the Grassmannian of complex $2$-planes of $\CC^m$ and 
five other exceptional spaces (see \cite{Be}, \emph{loc. cit.}) 
\end{ex}

        From the holonomy viewpoint, equivalent definitions are obtained as follows: 
A Riemannian manifold is quaternion \ka\  (resp. hyperk\"ahler) iff its holonomy 
is contained in $\mathrm{Sp}(n)\cdot \mathrm{Sp}(1)=\mathrm{Sp}(n)\times \mathrm{Sp}(1)/\ZZ_2$ (resp. $\mathrm{Sp}(n)$).

        On a quaternion  quaternionic Hermitian manifold, the usual \ka\ forms are only local: 
on any trivializing open set $U$, one has the $2$-forms $\omega_\al(\cdot,\cdot)=
g(I_a\cdot,\cdot)$. But the $4$-form 
$\omega=\sum_{\al=1}^3\omega_\al\wedge\omega_\al$ 
is global (because the transition functions of $H$ are in $\mathrm{SO}(3)$), nondegenerate   
and, if the manifold is quaternion \ka , parallel. Hence, it gives a nontrivial $4$-cohomology class, precisely 
$[\omega]=8\pi^2p_1(H)\in H^4(M,\RR)$ (\cite{Kr}). 

        To get a converse, let 
$\HA$ be the  algebraic ideal generated by $H$ in 
$\Lambda^2T^*M$ (by identifying, as usual, a local almost complex structure with  
the associated \ka\ $2$-form). It is a differential ideal if for any 
admissible basis of $H$, 
one has $d\omega_\al=\sum_{\be=1}^3\eta_{\al\be}\wedge\omega_\be$ for 
some local $1$-forms $\eta_{\al\be}$. Then we have:

\begin{te}\cite{Sw}~ \label{swa}
A quaternion Hermitian manifold of dimension at least $12$  
with closed $4$-form $\omega$ is quaternion \ka .

        A  quaternion Hermitian manifold of dimension $8$ is quaternion K\"a\-hler  
iff $\omega$ is closed and $\HA$ is a differential ideal.
\end{te}
Swann's proof uses representation theory. A more direct one can be found in 
\cite{AM1}. 
\begin{re}\label{conf}
It is important to note that the condition of being a differential ideal is 
conformally invariant and, moreover, invariant to different choices of 
admissible basis.
\end{re}

        For an almost almost quaternionic Hermitian manifold $(M,g,H)$, we         
define its structure tensor by
$$T^H=\frac{1}{12}\sum_{\al=1}^3[I_\al,I_\al].$$
Clearly, $T^H$ is zero if one can choose, locally, admissible basis formed by 
integrable almost complex structures. The \emph{O\-ba\-ta connection} $\nabla^H$ is 
then the unique connection which preserves $H$ and has torsion equal to $T^H$. 
It defines the fundamental $1$-form $\eta$ by the relation
$$\eta(x)=\frac{1}{(8(n+1)}\mathrm{trace}\{g^{-1}\nabla^H_Xg\}.$$
A direct (but lengthy) computation proves:
\begin{lm}\label{AM1} \cite{AM1} 
Let $(M,g,H)$ be a quaternion Hermitian manifold such that $\HA$ is a differential 
ideal. For any admissible basis of $H$, the following formulae for 
$T^H$ and $\nabla^H$ hold good:

\begin{equation}\label{t}
\begin{split}
T^H_XY=&\frac{1}{60}\sum_{\al=1}{3}\{[(5\f_\al+\rho_al)(I_\al X]I_\al Y-
[(5\f_\al+\rho_al)(I_\al Y]I_\al X+\\
+&4\omega_\al X,Y)g^{-1}(\rho_\al\circ I_\al)\},
\end{split}
\end{equation}
\begin{equation}\label{nab}
\begin{split}
(\nabla^H_Zg)(X,Y)=&\frac{1}{12}\{2\e(Z)g(X,Y)+\e(X)g(Y,Z)+\e(Y)g(X,Z)-\\
-&\sum_{\al=1}^3\e(I_\al X)g(Y,I_\al Z)-
\sum_{\al=1}^3\e(I_\al Y)g(X,I_\al Z)+\\
+&4\eta(Z)g(X,Y)\}
\end{split}
\end{equation}
where 
\begin{equation*}
\begin{split}
\f_\al=&2\eta_{[\be\ga]}\circ I_\al-\eta_{[\ga\al]}\circ I_\be-
\eta_{[\al\be]}\circ I_\ga,\\
\rho_\al=&-6\eta_{\al\al}+2\eta-3\eta_{(\al\be)}\circ I_\ga+
3\eta_{(\ga\al)}\circ I_\be\\
\e=&\sum_{\al=1}^3\eta_{[\al\be]}\circ I_\ga
\end{split}
\end{equation*}
the subscript $()$ (resp. $[]$) indicating symmetrization (resp. 
sqew-sym\-me\-tri\-za\-tion).
\end{lm}

        A most important geometric property of quaternion K\"ahler manifolds,  
partly motivating the actual interest in their study is:
\begin{te}\cite{Ber}
A quaternionic \ka\ manifold is Einstein.
\end{te}
We briefly sketch, following \cite{Be}, p. 403, S. Ishihara's proof. We fix a local 
admissible basis. Direct computations lead to the formulae:
\begin{equation*}
[R^g(X,Y),I_\al]=\sum_{\be=1}^3\eta_{\al\be}I_\be
\end{equation*}
with a sqew-symmetric matrix of $2$-forms $(\eta_{\al\be})$ which can be 
expressed in terms of  Ricci tensor as follows:
\begin{equation*}
\eta_{\al\be}(X,Y)=\frac{2}{n+2}Ric^g(I_\ga X,Y), \quad \dim M=4n.
\end{equation*}
From these one gets:

\begin{equation*}
\begin{split}
&g(R^g(X,I_1X)Z,I_1Z)+g(R^g(X,I_1X)I_2Z,I_3Z)+\\
+&g(R^g(I_2X,I_3X)Z,I_1Z)+g(R^g(I_2X,I_3X)I_2Z,I_3Z)=\\
=&\frac{4}{n+4}Ric^g(X,X)\norm{Z}^2=\frac{4}{n+4}Ric^g(Z,Z)\norm{X}^2
\end{split}
\end{equation*}
for any $X$ and $Z$, hence $Ric^g(X,X)=\la g(X,X)$ and $(M,g)$ is Einstein. 

        On the other hand, hyperk\"ahler manifolds have holonomy included in 
$\mathrm{Sp}(n)\subset \mathrm{SU}(2n)$, hence they are Ricci-flat, in particular Einstein 
(see \cite{Be}). 

        Although apparently hyperk\"ahler manifolds form a subclass of 
quaternion \ka\ ones, this is not quite true. Besides the holonomy argument, 
the following result motivates the dichotomy:
\begin{te} \cite{Ber} A quaternion K\"ahler manifold is Ricci-flat iff 
its reduced holonomy group is contained in $\mathrm{Sp}(n)$. And if it is not 
Ricci-flat, then it is de Rham irreducible.
\end{te}
From these results it is clear that when  
discussing  quaternion \ka\ manifolds, one is mainly interested in the 
non-zero scalar curvature.  

        Ricci flat quaternion K\"ahler manifolds are called 
\emph{locally hyperk\"ahler}. Similarly,  
P. Piccinni discussed in \cite{Pi2}, \cite{Pi} the class of \emph{locally 
quaternion K\"ahler manifolds}, having the \emph{reduced} holonomy group contained in 
$\mathrm{Sp}(n)\cdot \mathrm{Sp}(1)$ and proved:
\begin{pr}\cite{Pi2}
Any complete locally quaternion K\"ahler manifold with positive scalar 
curvature is compact, locally symmetric and admits a finite covering 
by a quaternion K\"ahler Wolf symmetric space.
\end{pr}

        As the local sections of $H$ are generally non-integrable, one 
cannot use the methods of complex geometry directly on quaternion \ka\ manifolds. 
However, one can construct an associated bundle whose total space is  
Hermitian. Let  $p:Z(M)\rightarrow M$ be the unit sphere subbundle of $H$. Its fibre 
$Z(M)_m$ is the set of all almost complex structures on $T_mM$. This is 
called the \emph{twistor bundle} of $M$. Using the Levi-Civita connection $\nabla^g$, 
one splits the tangent bundle of $Z(M)$ in horizontal and vertical parts. 
Then an almost complex structure $\JJ$ can be defined on $Z(M)$ as follows: each 
$z\in Z(M)$ represents a complex structure on $T_{p(z)}M$; as the horizontal  
subspace in $z$ is naturally identified with $T_{p(z)}M$, the action of $\JJ$ 
on horizontal vectors will be the tautological one, coinciding with the action of 
$z$. The vertical subspace in $z$ is isomorphic with the tangent space of the 
fibre $S^2$. Hence we let $\JJ$ act on vertical vectors as the canonical 
complex structure of $S^2$. Happily, $\JJ$ is integrable. Moreover:
\begin{te}\cite{Sa} Let $(M,g,H)$ be a quaternion \ka\ manifold with positive 
scalar curvature. Then $(Z(M),\JJ)$ admits a K\"a\-hler - Einstein metric with positive 
Ricci curvature with respect to which $p$ becomes a Riemannian submersion.
\end{te}

\section{Local and global $3$-Sasakian manifolds}\label{3ss}
We now describe the odd dimensional analogue, within the frame of contact  
geometry, of hyperk\"ahler manifolds, as well as a local version of it. 
We send the interested reader to  
the excellent recent survey \cite{BG2}, where also a rather exhaustive list 
of references is given and to \cite{OP2} for the local version.
\begin{de}
A $4n+3$ dimensional Riemannian manifold $(N,h)$ such that the c\^one metric 
$dr^2+r^2h$ on $\RR_+\times N$ is hyperk\"ahler is called a \emph{$3$-Sasakian}   
manifold.
\end{de}
This is equivalent to the existence of three 
mutually orthogonal unit Killing vector fields $\xi_1$, $\xi_2$, $\xi_3$, 
each one defining a
Sasakian structure (\emph{i.e.}: 
$\f_\al :=\nabla^h\xi_\al$ satisfies the differential equation 
$\nabla^h\f_\al=Id\otimes\xi_\al^\flat-h\otimes\xi_\al$)  
and related by: $$[\xi_1,\xi_2]=2\xi_3,
[\xi_2,\xi_3]=2\xi_1,
[\xi_3,\xi_1]=2\xi_2.$$ 
$3$-Sasakian manifolds are necessarily Einstein (\cite{Ka}) 
with positive scalar curvature and their Einstein
constant is $4n+2$.

        Starting with a $3$-Sasakian manifold $N$, one has to 
consider the foliation generated by the three structure vector fields $\xi_\al$. 
It is easy to compute the curvature of the leaves: it is precisely one. Hence, 
the leaves are spherical space forms. If the foliation is quasi-regular (it is enough 
to have compact leaves), 
then the quotient space is a quaternion \ka\ orbifold $M$ of positive 
sectional curvature (see \cite{BG3} for a thorough discussion about the 
geometry and topology of orbifolds and their applications in contact geometry). 
As all the geometric constructions we are interested in can be carried out in the 
category of orbifolds, one considers now the twistor space $Z(M)$. The triangle is   
closed by observing that, fixing one of the contact structures of $N$, one has 
an $S^1$-bundle $N\rightarrow Z(M)$  whose Chern class is, up to torsion, 
the one of an induced Hopf bundle (this is a particular case of a Boothby-Wang 
fibration, cf. \cite{Bl}).
Moreover, all three orbifold fibrations involved in this commutative triangle 
are Riemannian submersions.

        Conversely, given a positive quaternion \ka\ orbifold $(M,g,H)$, one 
constructs its \ka -Einstein twistor space (it will be an orbifold) and an  
$\mathrm{SO}(3)$-principal bundle over $M$. The total space $N$ will then be  a 
$3$-Sasakian orbifold which, as above, fibers in $S^1$ over $Z(M)$ closing 
the diagram. One of the deepest results in this theory was the 
determination of  conditions under which $N$ is indeed a manifold 
(cf. \cite{BGMR}). 

        A local version of $3$-Sasakian structure  
will be also useful in the sequel:
\begin{de}\cite{OP2}\label{3s} 
A Riemannian manifold $(N,h)$ is said to be a \emph{
locally $3$-Sasakian manifold} if a rank $3$ vector subbundle $\KK \subset TN$ is
given, locally spanned by an orthonormal triple $\xi_1,\xi_2,\xi_3$ of 
Killing 
vector fields satisfying:
\par (i)  $[\xi_\alpha,\xi_\beta]=2\xi_\gamma$ for $(\alpha , \beta , 
\gamma)=(1,2,3)$ and circular permutations.
\par (ii) Any two such triples $\xi_1,\xi_2,\xi_3$ and $\xi'_1,\xi'_2,\xi'_3$ 
are related 
on the intersections $U \cap U'$ of their definition open sets by matrices 
of functions with values in $\mathrm{SO}(3)$.
\par (iii) If $\f_\alpha = \nabla^h \xi_\alpha$ , 
 ($\alpha = 1,2,3$), then $ (\nabla^h_Y \f_\alpha)~Z=\xi_\alpha^\flat (Z)Y-
 h(Y,Z)\xi_\alpha,$
for any local vector fields $Y,Z$. 
\end{de}
Clearly, if $\KK$ can be globally trivialized with Killing vector fields as 
above, $(N,h)$ is $3$-Sasakian. It is easily seen that 
locally $3$-Sasakian manifolds share the local properties with the (global) 
$3$-Sasakian spaces: they are Einstein with positive scalar curvature; hence, 
by Myers' theorem we have
\begin{pr}
Complete locally and globally  $3$-Sasakian manifolds are compact. 
\end{pr}
But a specific property of the local case is:
\begin{pr}\cite{OP2}\label{flat}
The bundle $\KK$ of a locally $3$-Sasakian manifold is flat.
\end{pr}
\begin{proof}
Let $(\xi_1,\xi_2,\xi_3)$, $(\xi_1',\xi_2',\xi_3')$ be two local orthonormal 
triples of Killing fields trivializing $\KK$ on $U$, $U'$. Then, on 
$U\cap U'\neq \emptyset$ 
we have $\xi_\la'=f^\sigma_\la\xi_\sigma$. We shall show that $f^\sigma_\la$ are 
constant. Compute first the bracket
$$2\xi'_\nu=[\xi'\la,\xi'_\mu]=\{f^\rho_\la\xi_\rho(f_\mu^\sigma)-
f_\mu^\rho\xi_\rho(f_\la^\sigma)\}\xi_\sigma+
f^\rho_\la f_\mu^\sigma[\xi_\rho,\xi_\sigma].$$
From $(f^\mu_\la)\in \mathrm{SO}(3)$ and $[\xi_\rho,\xi_\sigma]=2\xi_\tau$ ($(\rho, \sigma, 
\tau) = (1,2,3)$ and cyclic permutations), we can derive:
$$f^\rho_\la f_\mu^\sigma[\xi_\rho,\xi_\sigma]=2\{f^\rho_\la f_\mu^\sigma-
f_\la^\sigma f_\mu^\rho\}\xi_\tau=2\xi'_\nu.$$
Hence  $$f^\rho_\la\xi_\rho(f_\mu^\sigma)-f_\mu^\rho\xi_\rho(f_\la^\sigma)=0.$$
Thus, for any $\la, \mu,\sigma=1,2,3$: $\xi'_\la(f^\mu_\sigma)-\xi'_\mu(f_\la^\sigma)=0$. 
It follows:
\begin{equation}\label{p}
\xi_\la(f_\sigma^\mu)-\xi_\mu(f^\la_\sigma)=0.
\end{equation}
Now we use the Killing condition applied to $\xi_\la'=f^\mu_\la\xi_\mu$: 
$$Y(f_\la^\mu)h(\xi_\mu,Z)+Z(f_\la^\mu)h(\xi_mu,Y)=0, \quad Y,Z\in \mathcal{X}(M)$$
which yields, on one hand $Z(f_\la^\mu)=0$ for any $Z\perp span\{\xi_1,\xi_2,\xi_3\}$ 
and, on the other hand
\begin{equation}
\xi_\rho(f_\la^\sigma)+\xi_\sigma(f_\la^\rho)=0.
\end{equation}
This and \eqref{p} imply $\xi_\sigma(f_\la^\rho)=0$ and the proof is complete.

\end{proof}
The vector bundle $\KK$ generates a $3$-dimensional foliation that, for 
simplicity, we equally denote $\KK$. It can be shown that $\KK$ is Riemannian. 
As in the global case, if the leaves of $\KK$ are compact, the leaf space 
$M=N/\KK$ is a compact orbifold. The metric $h$ projects to a metric $g$ on 
$P$ making 
the natural projection $\pi$ a Riemannian submersion with totally geodesic 
fibers. The locally 
defined endomorphisms $\f_\la$ can be projected on $M$ producing locally 
defined almost complex structures: $J_\al X_{\pi(x)}=\pi_*(\f_\al(\tilde X_x))$, 
where $\tilde X$ is the horizontal lift of $X$ w.r.t. the submersion. 
As $\f_\al\circ\f_\be=-\f_\ga+\xi_\al\otimes\xi_\be^\flat$, $P$ can be 
covered with open sets endowed with local almost hyperhermitian structures 
$\{J_\al\}$. As the transition functions of $\KK$ are in $\mathrm{SO}(3)$, so are 
the transition functions of the bundle $\mathcal{F}$ locally spanned 
by the $\f_\al$. Hence, two different almost hyperhermitian structures are 
related on their common domain by transition functions in $\mathrm{SO}(3)$. This means that 
the bundle $H$ they generate is quaternionic. Using the O'Neill formulae, it 
is now  seen, as in the global case, that $(M,g,H)$ is a quaternion 
K\"ahler orbifold. Summing up we can state:
\begin{pr}\cite{OP2}\label{fib}
Let $(N,h,\KK)$ be a locally $3$-Sasakian manifold such that $\KK$ has compact 
leaves. Then the leaf space $M=N/\KK$ is a quaternion K\"ahler orbifold with 
positive scalar curvature and the natural projection $\pi:N\rightarrow M$ is 
a Riemannian, totally geodesic  submersion which fibers are 
(generally inhomogeneous) $3$-dimensional spherical space forms.
\end{pr}
\begin{re}
P. Piccinni proved in \cite{Pi2} that some global $3$-Sasakian manifolds also 
project over local quaternion K\"ahler manifolds with positive scalar curvature. 
\end{re}
A further study of the (supposed compact) leaves of $\KK$ will show a very specific property of 
locally $3$-Sasakian manifolds. To this end, we recall, following \cite{Sc}, 
some aspects of the classification of $3$-dimensional spherical space forms $S^3/G$, with 
$G$ a finite group of isometries of $S^3$, hence  
a finite subgroup of $\mathrm{SO}(4)$. The finite subgroups of $S^3$ are known: they are 
cyclic groups of any order or binary dihedral, tetrahedral, octahedral, icosahedral and, 
of course, the identity. In all these cases, $S^3/G$ is a homogeneous 
$3$-dimensional space form carrying an induced (global) $3$-Sasakian structure, 
see \cite{BGM0}. The other finite subgroups of $\mathrm{SO}(4)$, not contained in 
but acting freely on $S^3$, are characterized by being conjugated in 
$\mathrm{SO}(4)$ to a subgroup of $\Gamma_1=\mathrm{U}(1)\cdot \mathrm{Sp}(1)$ or 
$\Gamma_2=\mathrm{Sp}(1)\cdot \mathrm{U}(1)$. Observe that the right (resp. left) isomorphism 
between $\HH$ and $\CC^2$ induces an isomorphism between  $\Gamma_1$ (resp. 
$\Gamma_2$)  and $\mathrm{U}(2)$. Hence, any finite subgroup $\Gamma$ of $\Gamma_1$ or $\Gamma_2$ will 
preserve two structures of $S^3$: the locally $3$-Sasakian structure 
induced by the hyperhermitian structure of $\CC^2$ and a global 
Sasakian structure induced by some complex Hermitian structure of 
$\CC^2$ belonging to the given hyperhermitian one. Moreover, altering 
$\Gamma$ by conjugation in $\mathrm{SO}(4)$ does not affect the above preserved structures; 
only the global Sasakian structure will come from a hermitian structure of $\RR^4$ 
conjugate with the standard one. Altogether, we obtain:

\begin{pr}\cite{OP2}\label{gog}
On any locally $3$-Sasakian manifold, the compact leaves of $\KK$ are 
locally $3$-Sasakian $3$-dimensional space-forms carrying a global  
almost Sasakian structure.
\end{pr}
We end with another consequence of Proposition \ref{gog}:
\begin{co} \cite{OP2}\label{pul}
Let $\tilde \KK\rightarrow \tilde N$ be the pull-back of the bundle 
$\KK\rightarrow N$ to the universal Riemannian covering space of a locally $3$-Sasakian manifold. 
Then $\tilde \KK$ is globally trivialized by a global $3$-Sasakian structure on 
$\tilde N$.
\end{co}
\begin{proof}
By Proposition \ref{flat}, the bundle $\tilde\KK\rightarrow \tilde N$ is 
trivial. However, this is not enough 
to deduce that the trivialization can be realized with Killing fields 
generating a $\mathrm{su}(2)$ algebra. \emph{E.g.} the inhomogeneous 
$3$-di\-men\-si\-o\-nal 
spherical space forms are parallelizable but locally, not globally 
$3$-Sasakian. To overcome this difficulty, start with the induced 
locally $3$-Sasakian structure of $\tilde N$. Let $X_1$ be the global 
Sasakian structure of $\tilde N$ provided by Proposition \ref{gog} and 
consider an open set $\tilde U$ on which $\tilde\KK$ is trivialized by 
a local $3$-Sasakian structure incuding $X_1$. 

        The manifold $\tilde N$ is 
simply connected and Einstein, hence analytic (see \cite{Be}, Theorem 5.26). 
By a result of Nomizu 
(cf. \cite{No}) each local Killing vector field on $\tilde N$ can be 
extended uniquely to the whole $\tilde N$. We thus extend the above 
three local Killing fields. Clearly, the extension $Y_1$ of $X_1$ 
coincides with $X_1$.  
The extension $Y_2$ of ${X}_2$ is thus orthogonal to $Y_1$ and belongs to 
$\tilde \KK$ in every point of $\tilde N$. It follows from Proposition \ref{flat} 
that ${Y}_2$ is a global Sasakian structure. Now $Y_3=\frac{1}{2}[Y_1,Y_2]$ 
completes the desired global $3$-Sasakian structure.
\end{proof}

\section{Quaternion Hermitian Weyl and hyperhermitian Weyl manifolds}

        We now arrive to the structures giving the title of this survey. We consider 
$4n$-dimensional quaternion Hermitian manifolds and let the metric vary in 
its conformal class. In this setting, the natural connection to work with is no more 
the Levi-Civita connection, but a Weyl connection which has to be compatible 
with the quaternionic structure too.

\subsection{Definitions. First properties}

        Let  $(M^{4n},c,H)$, $n \geq 2$ be a  conformal   manifold 
endowed with a quaternionic bundle $H$ such that $(M,g,H)$ is quaternion Hermitian 
for each $g\in c$. 

\begin{de} $(M^{4n},H,c,D)$ is said 
\emph{ quaternion-Hermitian-Weyl} if:  

        1)~ $(M,c,D)$ is a Weyl manifold; 

        2)~ $(M,g,H)$ is quaternion-Hermitian for any $g\in c$;

        3)~ $DH=0$, \emph{i.e.} $DI_\al=a^\be_\al \otimes I_\be$ with sqew-symmetric 
matrix of one-forms $(a_\al^\be)$ for any admissible basis of $H$.

        $(M^{4n},c,H,D)$ is said \emph{hyperhermitian Weyl} if it satisfies condition 1)  
and:        

        2')~ $(M,g,H)$ is hyperhermitian for any $g\in c$;

        3')~ $DI=0$ for any section of $H$.

\end{de}

The above definition is clearly inspired by the complex case, where the 
theory of 
Hermitian-Weyl (locally conformal \ka ian in other terminology) is 
widely studied (see \cite{DO} for a recent survey).  
Indeed, the following equivalent definition is available:
\begin{pr}\cite{PPS}
$(M^{4n},c,H,D)$ is quaternion-Hermitian-Weyl\\ (resp. hyperhermitian Weyl) 
if and only if $(M,g,H)$ is 
locally conformally quaternion \ka\ (resp. locally conformally hyperk\"ahler)  
(i.e. $g_{\vert_{U_i}} = e^{f_i}g'_i$, where the $g'_i$ are quaternion 
K\"ahler (resp. hyperk\"ahler) over open neighbourhoods $\{U_i\}$
covering $M$) for each $g\in c$.
\end{pr}
\begin{proof}
Let  $(M^{4n},c,H,D)$ be quaternion-Hermitian Weyl. Fix a metric $g\in c$ 
and choose an open set $U$ on which 
$H$ is trivialized by an admissible basis $I_1,I_2,I_3$. Then  $Dg=\theta_g\otimes g$ 
together with condition 2) of the definition imply $D\omega_\al=
\theta\otimes\omega_\al+a_\al^\be\otimes \omega_\be$, hence 
\begin{equation}\label{cc}
d\omega_\al= \theta\wedge\omega_\al+a_\al^\be\wedge \omega_\be. 
\end{equation}
This implies 
that $\HA$ is a differential ideal 
and, on the other hand, the derivative of the fundamental four-form is 
$d\omega=\theta_g\wedge\omega$. Differentiating here we get 
$0=d^2\omega=d\theta_g\wedge\omega$. As $\omega$ is nondegenerate, this means 
$d\theta_g=0$. Consequently, locally, on some open sets $U_i$, 
$\theta_g=df_i$ for some differentiable functions defined on $U_i$. It is now easy to 
see that for each $g'_i=e^{-f_i}g_{\vert_{U_i}}$, the associated $4$-form  
is closed, hence, taking into account Proposition \ref{swa} and Remark 
\ref{conf}, the local metrics $g'_i$ are quaternion \ka . 

        Conversely, starting with the local quaternion \ka\ metrics 
$ g'_i=e^{-f_i}g_{\vert_{U_i}}$, define $(\theta_g)_{\vert_{U_i}}=df_i$. 
It can be seen that these local one forms glue together to a global, closed 
one-form and $d\omega=\theta_g\wedge\omega$. Then construct the Weyl connection 
associated to $g$ and $\theta_g$:
\begin{equation}\label{cw}
D=\nabla^g-{\frac{1}{2}}\{\theta\otimes Id+Id\otimes\theta-
g\otimes\theta_g^{\sharp}\}.
\end{equation}
A straightforward computation shows that $D$ has the requested properties.

        The proofs for the global case are completely similar.
\end{proof}
\begin{co}\label{cor}
A quaternion Hermitian manifold $(M,g,H)$ is quaternion Hermitian Weyl if and 
only if there exist a $1$-form $\theta$ (necessarily closed) such that 
the fundamental $4$-form $\omega$ satisfy the integrability condition 
$d\omega=\theta\wedge\omega$. In particular, $(M,g,H)$ is 
quaternion K\"ahler if and only if $\theta=0$.
\end{co}
The form $\theta$ is the Higgs field associated to the Weyl manifold 
$(M,c,D)$. But in this context, we shall prefer to call it 
\emph{the Lee form} (see \cite{DO} for a motivation). 

        As on a simply connected manifold any closed form is exact we derive:
\begin{co}
A quaternion Hermitian Weyl ( hyperhermitian Weyl) manifold which is not 
globally conformal quaternion K\"ahler ( hy\-per\-k\"ah\-ler) cannot be simply 
connected. 

        The universal Riemannian covering space of a 
quaternion Hermitian Weyl ( hyperhermitian Weyl) manifold is 
globally conformal quaternion K\"ahler ( hyperk\"ahler).
\end{co}

\begin{ex}\label{ex}
We give here just one example of compact hyperhermitian Weyl manifold and 
leave the description of other examples for the end of the paper, 
following the structure of quaternion Hermitian Weyl and hyperhermitian 
Weyl manifolds. 

        The standard quaternionic Hopf manifold is $H^n_\HH=\HH-\{0\}/\Gamma_2$, where 
$\Gamma_2$ is the cyclic group generated by the quaternionic automorphism 
$(q_1,...,q_n)\mapsto (2q_1,...,2q_n)$. The hypercomplex structure 
of $\HH^n$ is easily seen to descend to $H_\HH^n$. Moreover, 
the globally conformal quaternion K\"ahler metric $(\sum_iq_i\ov{q}_i)^{-1}
\sum_idq_i\otimes d\ov{q}_i$ on $\HH^n-\{0\}$ is invariant to the action of 
$\Gamma_2$, hence induces a locally conformally hyperk\"ahler metric on 
the Hopf manifold. Note that, as in the 
complex case, $H^n_\HH$ is diffeomorphic with a product of spheres 
$S^1\times S^{4n-1}$. Consequently, 
its first Betti number is $1$ and it cannot accept any hyperk\"ahler metric.
\end{ex}
Before going over, let us note the following result:
\begin{pr}\cite{OP2}
A quaternion Hermitian manifold $(M,g,H)$ admits a unique 
quaternion Hermitian Weyl structure.
\end{pr}
\begin{proof}
We have to prove that there exists a unique torsion free connection preserving 
both $H$ and $[g]$. Indeed, if $D_1$, $D_2$ are such, let $\theta_1$, $\theta_2$ 
be the associated Lee forms. Then the fundamental $4$-form $\omega$ satisfies 
\begin{equation}\label{inj}
d\omega=\theta_1\wedge\omega=\theta_2\wedge\omega. 
\end{equation}
Using the operator $L:\Lambda^1T^*M
\rightarrow \Lambda^5T^*M$, $L\al=\al\wedge\omega$, \eqref{inj} yields   
$L(\theta_1-\theta_2)=0$. But $L$ is injective, because it is related to 
its formal adjoint $\Lambda$ by 
$\Lambda L=(n-1)Id$. Hence $\theta_1=\theta_2$.  Finally, formula 
\eqref{cw} proves that $D_1=D_2$.
\end{proof}
\begin{re}\cite{Pi2}
For hyperhermitian Weyl manifolds, this uniqueness property is implied  
by  the 
characterization of the Obata connection as the unique 
torsion-free hypercomplex connection. It must then coincide with our 
Weyl connection $D$. In general, the set of torsion-free quaternionic 
connections has an affine structure modelled on the space of $1$-forms. 
However, only one torsion-free connection can preserve a given conformal 
class of hyperhermitian metrics. This follows from the fact that 
the exterior multiplication with the fundamental four-form  of the 
metric maps injectively $\Lambda ^1(T^*M)$ into $\Lambda^5(T^*M)$. 
\end{re}

Note that the connection $D\vert_{U_i}$ is in fact the Levi Civita 
connection of the local quaternion K\"ahler metric $g'_i$. As quaternion-K\"ahler 
manifolds are Einstein, we obtain the following fundamental result: 
\begin{pr}\cite{PPS}
Quaternion Hermitian Weyl manifolds are Einstein Weyl.
\end{pr}
Hence, as $d\theta=0$, \emph{i.e.} the Weyl structure $(M,c,D)$ is 
closed and not exact, because the $D$ is the Levi-Civita connection of 
\emph{local} metrics (the Weyl structure is only locally exact), 
the quoted Theorem \ref{gg} of P. Gauduchon implies:
\begin{pr}\label{par}\cite{PPS}
On any compact quaternion-Weyl (hyperhermitian Weyl) manifold which is not globally 
conformal quaternion K\"ahler (hyperk\"ahler) there exists a representative 
$g\in c$ (the Gauduchon metric) such that the associated Lee form 
$\theta_g$ be $\nabla^g$-parallel.
\end{pr}
In the sequel, the parallel Lee form of the Gauduchon metric will always be 
supposed of unit length.
\begin{co}\label{va}
Let $g$ be the above metric with parallel Lee form on a compact 
hyperhermitian Weyl manifold and $\{I_\al\}$ an adapted hyperhermitian 
structure. Then $(g,I_\al)$ are Vaisman structures on $M$ (cf. \cite{DO}).
\end{co}
\begin{pr}\cite{OP1}
On a compact quaternion Weyl manifold which is not globally 
conformal quaternion K\"ahler, the local quaternion K\"ahler metrics $g'_i$ are 
Ricci-flat.
\end{pr}
\begin{proof}
This result follows directly from Theorem \ref{gg}, 2), but we prefer to 
give here a direct proof, adapted to our situation. 

        On each $U_i$, the relation between the scalar curvatures  
$\mathrm{Scal}'_i$ and $Scal$ 
of $g'_i$ and $g$ is (cf. \cite{Be}, p. 59):
$$\mathrm{Scal}'_i=e^{-f_i}\left\{\mathrm{Scal}_{|U_i}-\frac{(4n-1)(2n-1)}{2}\right\}.$$
Hence $\mathrm{Scal}'_i$ is constant. If $\mathrm{Scal}'_i$ is not identically zero, 
differentiation of the above identity yields:
$$\theta_{|U_i}=d\log \left\{\mathrm{Scal}_{|U_i}-\frac{(4n-1)(2n-1)}{2}\right\}.$$ 
As both $\theta$ and $\mathrm{Scal}$ are global objects on $M$, it follows that  
$\theta$ is exact, contradiction. But if $\mathrm{Scal}'_i=0$ on some $U_i$, then 
$Scal=\mathrm{Scal}_{|U_i}=\frac{(4n-1)(2n-1)}{2}$, constant on $M$. This proves that 
$\mathrm{Scal}'_i=0$ on each $U_i$.
\end{proof}
\begin{re}\label{int} 

The above result says that quaternion-Her\-mi\-ti\-an Weyl manifolds are 
locally conformally locally hyperk\"ahler. In particular, the open subsets 
$U_i$ can always be taken simply connected and endowed with admissible basis 
made of by integrable, parallel almost complex structures. But this does not 
mean that $M$ would be a locally conformal K\"ahler manifold, because a 
global K\"ahler structure might not exist.
\end{re}

Another characterization, using the differential ideal $\HA$ is the following 
(recall that the differential ideal condition is conformally invariant, so 
one can speak about the differential ideal of a conformal manifold):
\begin{te}\cite{ABM}\label{laba}
A quaternionic conformal manifold $(M,c,H)$ of dimension at least $12$ 
is quaternion-Hermitian Weyl 
if and only if $\HA$ is a differential ideal.
\end{te}
The following result is essential in the author's proof, also motivating 
the restriction on the dimension:
\begin{lm}\cite{ABM}
Let $(M,g,H)$ be an almost hyperhermitian manifold with $\dim M\geq 12$. 
Suppose 
$\sum_{\al=1}^3\phi_\al\wedge\omega_\al=0$ 
for some $2$-forms $\phi_\al$. Then there exists the sqew-symmetric matrix  
of real functions $f_{\al\rho}$ such that 
$\phi_\al=\sum_{\rho\neq\al}f_{\al\rho}\omega_\rho$.
\end{lm}
\begin{proof}
Let $F_\al$ be the $1-1$ tensor fields metrically equivalent with the 
$2$-forms $\phi_\al$. The identity in the statement can be rewritten as:
\begin{equation}\label{iu}
\begin{split}
\sum_{\rho=1}^3\{&-\phi_\rho(X,Y)I_\rho Z+\phi_\rho(X,Z)I_\rho Y +
\omega_\rho(Y,Z)F_\rho X-\\
&-\phi_\rho(Y,Z)I_\rho X+ \omega_\rho(Z,X)F_\rho Y+
\omega_\rho(X,Y)F_\rho Z\}=0.
\end{split}
\end{equation}
Let now $X$ be unitary, fixed. In the orthogonal complement of 
$\HH X=\{X, I_1X, I_2X,I_3X\}$ we choose a unitary $Z$ and let $Y=I_\al Z$. 
With these choices, the above identity reads: 
$$F_\al(X)= \sum_{\rho=1}^3\{\phi_\rho(X,I_\rho Z)I_\rho Z -
\phi_\rho(X,Z)I_\rho I_\al Z\}+\sum_{\rho=1}^3 \phi_\rho(I_\rho Z,Z)I_\rho X.$$
Here we use the assumption $n\geq 3$ to obtain:
$$F_\al(X)=\sum_{\rho=1}^3 \phi_\rho(I_\rho Z,Z)I_\rho X,$$
hence $\phi_\al$ have the form 
$\phi_\al=\sum_{\rho\neq\al}f_{\al\rho}\omega_\rho$ which, introduced in the 
equation \eqref{iu}, gives: 
\begin{equation*}
\begin{split}
\sum_{\al=1}^3f_{\al\al}\{-\omega_\al(X,Y)I_\al Z+\omega_\al(X,Z)I_\al Y-
\omega_\al(Y,Z)I_\al X\}+\\
+\sum_{\rho\neq\al}(f_{\al\rho}+f_{\rho\al})
\{\omega_\al(X,Y)I_\rho Z+\omega_\al(X,Z)I_\rho Y+
\omega_\al(Y,Z)I_\rho X\}=0.
\end{split}
\end{equation*}
Again using $n\geq 3$, we may choose $Y$ and $Z$ orthogonal to $\HH X$ and get:
$$-f_{\al\al}\omega_\al(Y,Z)- 
\sum_{\rho\neq\al}(f_{\al\rho}+f_{\rho\al})\omega_\al(Y,Z)=0.$$
Now it remains to take $Z=I_\al Y$ to derive the sqew-symmetry of $(f_{\al\rho})$. 
\end{proof}

\begin{proof} (of Theorem \ref{laba}). 
Fix $g\in c$ and an admissible basis for $H$. Starting from  equations 
\eqref{t}, \eqref{nab} and 
$d\omega_\al=\sum_{\be=1}^3\eta_{\al\be}\wedge\omega_\be$, one can derive 
the following formula:
\begin{equation}\label{cu}
d\omega_\al=\eta_\ga\wedge\omega_\be-\omega_\be\wedge\omega_\ga+\frac{1}{3}
\eta\wedge\omega_al,
\end{equation}
where $2\eta:=\eta_{\be\ga}-\eta_{\ga\be}$. 
After differentiating  \eqref{cu} we get:
$$\frac{1}{3}d\eta\wedge\omega_\al+(d\eta_\ga+\eta_\al\wedge\eta_\be)\wedge\omega_\be-
(d\eta_\be+\eta_\ga\wedge\eta_\al)\wedge\omega_\ga=0.$$
The previous Lemma applies and provides:
\begin{equation*}
\begin{split}
\frac{1}{3}\eta=&f_{\al\be}\omega_\be+f_{\be\ga}\omega_\ga\\
d\eta_\ga+\eta_\al\wedge\eta_\be=&f_{\be\al}\omega_\al+f_{\be\ga}\omega_\ga\\
-d\eta_\be+\eta_\ga\wedge\eta_\al=&f_{\ga\al}\omega_\al+f_{\ga\be}\omega_\be
\end{split}
\end{equation*}
This yields $d\eta=0$ and $d\eta_\al+\eta_\be\wedge\eta_\ga=f\omega_\al$ with 
$f$ not depending on $\al$. Hence, locally $f=d\sigma$ and we have 
$$d\omega_\al=\eta_\ga\wedge\omega_\be-\eta_\be\wedge\omega_\ga-
\frac{1}{3}d\sigma\wedge\omega_\al,$$
an equation similar to \eqref{cc}. The rest and the converse are obvious.
\end{proof}
\begin{re}
It is still unknown if this result is true in dimension $8$ too.
\end{re}
\begin{re}
For quaternion Hermitian manifolds, various adapted 
canonical connections were  
introduced by V. O\-pro\-iu, M. Obata and others. A unified treatement 
can be found in some recent papers of D. Alekseevski, E. Bonan, S. Marchiafava 
(see \emph{e.g.} \cite{AM2} and the references therein). 
In particular, in \cite{Ma3}, one finds  
a characterization of hyperhermitian Weyl manifolds in terms of canonical 
connections and structure tensors of the subordinated quaternionic Hermitian 
structure.
\end{re}
        We end this section with a characterizations of quaternion K\"ahler 
manifolds among (non compact) quaternion Hermitian Weyl manifolds by means of 
submanifolds (compare with \cite{V2} for the complex case):
\begin{pr}\cite{OP1}
A quaternion Hermitian Weyl manifold $(M,g,$ $H)$ of dimension at least $8$ 
is quaternion K\"ahler if and only if through 
each point of it passes a totally geodesic submanifold of real dimension $4h
\geq 8$ which is quaternion K\"ahler with respect to the structure 
induced by $(g,H)$.
\end{pr}
\begin{proof}
On a given submanifold of $M$, locally one can induce the metric $g$ and the 
quaternion K\"ahler one $g_i'$. Correspondingly, there are two second fundamental 
forms $b$ and $b_i'$. As $g$ and $g_i'$ are conformally related on $U_i$, the relation 
between $b$ and $b_i'$ is  
$$b'_i=b+\frac{1}{2}g\otimes T^\nu,$$ 
where $T^\nu$ is the part of $T$ normal to the submanifold. Now let 
$x\in M$ and let $j:Q\rightarrow M$ be a quaternion K\"ahler submanifold through $x$ 
as stated. We have $j^*d\omega=0$. From  $d\omega=\theta\wedge\omega$ we then 
derive $j^*\theta\wedge j^*\omega=0$. But rank $j^*\omega=4h\geq 8$, hence 
$j^*\theta=0$ meaning that $T$ is normal to $Q$: $T=T^\nu$. On the other hand, 
the same relation $j^*\theta=0$ shows that $Q\cap U_i$ is a quaternion 
K\"ahler submanifold of  the quaternion K\"ahler manifold 
$(U_i,H_{|U_i},g'_i)$. As quaternion submanifolds of quaternion K\"ahler 
manifolds are totally geodesic, $Q\cap U_i$ is totally geodesic in $U_i$ 
with respect to $g'_i$. It follows $2b=-g\otimes T$ on $Q\cap U_i$. But 
$b$ is zero from the assumption ($Q$ is totally geodesic with respect to $g$). 
This yields $T=0$ on $Q\cap U_i$, in particular $T_x=0$. Since $x$ was arbitrary in $M$,   
$T=0$ on $M$ proving that $(M,g,H)$ is quaternion K\"ahler.

        For the converse, just take $Q=M$.
\end{proof}

We end this general presentation with a recent result which makes quaternion Hermitian Weyl manifolds interesting for physics. We first recall (sending to \cite{GP} and \cite{Iv} for details and further references)  
the notion of \emph{quaternionic K\"ahler} (resp. \emph{hyperk\"ahler}) \emph{manifold with torsion}, briefly QKT (resp. HKT) manifolds. Let  $(M,g,H)$ be a quaternionic Hermitian (resp. hyperhermitian) manifold. It is called QKT (resp. HKT) manifold if it admits a metric quaternionic (resp. hypercomplex) connection $\nabla$ with totally skew symmetric torsion tensor which is, moreover, of type $(1,2)+(2,1)$ w.r.t. each local section $I_\al$, that is it satisfies:
$$T(X,Y,Z)=T(I_\al X,I_\al Y,Z)+T(I_\al X,Y,I_\al Z)+T(X,I_\al Y,I_\al Z),$$
where $T(X,Y,Z)=g(Tor^\nabla(X,Y),Z)$ and $Tor^\nabla(X,Y)=\nabla_XY-\nabla_YX-[X,Y]$. The holonomy of such a connection is contained in $\mathrm{Sp}(n)\cdot \mathrm{Sp}(1)$. These structures appear naturally on the target space of $(4,0)$ supersymmetric two-dimensional sigma model with Wess-Zumino term and seem to be of growing interest for physicists. Let us introduce the $1$-forms:
$$t_\al(X)=-\frac 12\sum_{i=1}^{4n}T(X,e_i,I_\al e_i),\quad \al=1,2,3.$$
Then  the $1$-form $t=I_\al t_\al$ is independent on the choice of $I_\al$.  
We can now state:
\begin{pr}\cite{Iv}
Every quaternion Hermitian Weyl (resp. hyperhermitian Weyl) manifold admits a QKT (resp. HKT) structure.

Conversely, a $4n$ dimensional ($n>1$) QKT manifold $(M,g,H,\nabla)$ is quaternion Hermitian Weyl if and only if:
$$T=\frac{1}{2n+1}\sum_\al t_\al\wedge\omega_\al \;\; \text{and} \;\; dt=0.$$
\end{pr}

\subsection{The canonical foliations}

        From now on $(M,c,H,D)$ will be  compact, non globally conformal 
quaternion K\"ahler. According to Proposition 
\ref{par}, we let $g\in c$ be the Gauduchon metric whose Lee form 
$\theta:=\theta_g$ 
is parallel w.r.t. the Levi-Civita connection $\nabla:=\nabla^g$. Hence 
we look at the quaternion Hermitian  manifold $(M,g,H)$. We also suppose 
$\theta\neq 0$ meaning that $M$ is  not quaternion 
K\"ahler, see Corollary \ref{cor}. We recall that, being parallel, 
we can suppose $\theta$ 
normalised, \emph{i.e.} $\mid\theta\mid=1$.  We denote $T:=\theta^\sharp$  
and let $T_\al=I_\al T$ and $\theta_\al=\theta\circ I_\al$.

        The following proposition gathers the computational formulae we need:

\begin{pr}\cite{OP1}
Let $(M,g,H)$ be a compact quaternion Hermitian Weyl manifold and $\{I_1,I_2, 
I_2\}$ a local admissible basis of $H$ with $I_\al$ integrable  
and parallel (as in remark \ref{int}). The following formulae hold good:
\begin{gather}
\label{a}\Ll_TI_\al=0,\quad \Ll_Tg=0,\quad \Ll_T\omega_\al=0, 
        \quad \Ll_T\omega=0\\
\label{b}\nabla I_\al=\frac{1}{2}\{Id\otimes\theta_\al-I_\al\otimes\theta-
        \omega_\al\otimes T+g\otimes  T_\al\}\\
\label{c}\Ll_{T_\al}I_\al=0, \quad \Ll_{ T_\al}I_\be=I_\ga, 
        \quad \Ll_{ T_\al}g=0\\
\label{d}[T, T_\al]=0, \quad [T_\al,  T_\be]= T_\ga \\
\label{e}\nabla \theta_\al=\frac{1}{2}\{\theta\otimes\theta_\al-\theta_\al\otimes\theta- 
        \omega_\al\}\\
\label{f}d\theta_\al=-\omega_\al+\theta\wedge\theta_\al\\
\label{g}\Ll_{T_\al}\omega_\al=0, \quad \Ll_{T_\al}\omega_\be=\omega_\be, \quad 
        \Ll_{T_\al}\omega=0
\end{gather}
where $\Ll$ is the operator of Lie derivative.
\end{pr}
The proof is by direct computation and mimics the corresponding one for Vaisman 
manifolds, see \cite{DO}. In particular, from \eqref{a} and \eqref{c},  
we obtain according to \cite{P1}:
\begin{co}
The vector fields $T$ and $T_\al$ are infinitesimal automorphisms of the 
quaternion Hermitian structure. 
\end{co}

There are two interesting foliations on any compact 
quaternion Hermitian Weyl manifold: 
\begin{itemize}
\item the $(4n-1)$-dimensional $\F$, spanned by the kernel of $\theta$ and 
\item the $4$-dimensional $\D$, locally generated by $T, T_1,T_2,T_3$. 
\end{itemize}
Here are their properties:
\begin{pr}\cite{OP2},\label{ff}
On a compact quaternion Hermitian Weyl manifold, 
$\F$ is a Riemannian, totally geodesic foliation. Its leaves have an 
induced locally $3$-Sasakian structure.
\end{pr}
\begin{proof}
The first statement is a consequence of \eqref{a}. As for the second one, 
the bundle $\KK$ is locally generated by the (rescaled to be unitary) 
local vector fields $T_\al$. Indeed, they are Killing by the last equation 
of \eqref{c}; the first condition of definition \ref{3s} is given by \eqref{d}; 
the transition functions of $\KK$ are in $\mathrm{SO}(3)$ because the transition 
functions of $H$ are so; finally, condition 3) of the definition is 
implied by \eqref{e}. 
\end{proof}
\begin{co}\cite{OP1}
On a compact hyperhermitian Weyl manifold, $\F$ is a Riemannian, totally geodesic 
foliation whose leaves have an induced (global) $3$-Sasakian structure.
\end{co}
\begin{pr}\cite{OP1}\label{qq}
On a compact quaternion Hermitian Weyl manifold, the foliation $\D$ is 
Riemannian, totally geodesic. Its leaves are conformally flat $4$-manifolds 
$(\HH-\{0\})/G$, with $G$ a discrete subgroup of $\mathrm{GL}(1,\HH)\cdot \mathrm{Sp}(1)$ 
inducing an integrable (in the sense of G-structures) quaternionic structure.
\end{pr}
\begin{proof}
Let $X$  be a leaf of $\D$ and let the superscript $'$ refer to 
restrictions of objects from $M$ to $X$. A local 
orthonormal basis of tangent vectors for $X$ is provided by $\{T',T_1',
T_2',T_3'\}$. As $X$ is totally geodesic, $\nabla'\theta'=0$ and   
a direct computation of the curvature tensor of the Weyl connection $R^D$ 
on this basis proves $R^D=0$ on $X$. Hence $X$ is conformally flat and 
the curvature tensor of the Levi-Civita connection is
\begin{equation}\label{cur}
\begin{split}
R'(U,Y)Z=&\theta'(U)\theta'(Z)Y-\theta'(Y)\theta'(Z)U-\theta'(U)g'(Y,Z)T'+\\
+&\theta'(Y)g'(U,Z)T'+g'(Y,Z)U-g'(U,Z)Y.
\end{split}
\end{equation}
It follows that the Ricci tensor $Ric'=g'-\theta'\otimes\theta'$ is 
$g'$-parallel and, on the other hand, the sectional curvature is non-negative and 
strictly positive on any plane of the form $\{T'_\al,T'_\be\}$. Now recall that  
the universal Riemannian covering spaces of conformally flat Riemannian 
manifolds with parallel Ricci tensor were classified in \cite{L}. By the above 
discussion and the reducibility of $X$ (due to $\nabla'T'=0$), the only class 
fitting from Lafontaine's classification is that with universal cover 
$\RR^4-\{0\}$ equipped with the conformally flat metric written 
in quaternionic coordinate $(h\ov{h})^{-1}dh\otimes d\ov{h}$. We still have to 
determine the allowed deck groups.

        Happily,  Riemannian manifolds with such universal cover 
were studied in \cite{G1} and, in arbitrary dimension, in 
\cite{RV}. Here it is proved that equation \eqref{cur} forces the 
deck group of the covering to contain only conformal transformations of 
the form (in real coordinates)
$\tilde{x}^i=\rho a^i_jx^j$ where $\rho>0$ and $(a_j^i)\in \mathrm{SO}(4)$. This 
leads to the following form of $G$: 
\begin{equation}\label{54}
G=\{ht_0^k\; ;\; h\in G_0, k\in \ZZ\}
\end{equation}
where $t_0$ is a conformal transformation of maximal module $0<\rho<1$ and 
$G'$ is one of the finite subgroups of $\mathrm{U}(2)$ listed in \cite{Kat}. Finally, 
as $\mathrm{CO}^+(4)\simeq \mathrm{GL}(1,\HH)\cdot \mathrm{Sp}(1)$, $X$ has an induced integrable 
quaternionic structure.
\end{proof}
\begin{co}\cite{OP1} \label{qqq}
On a compact hyperhermitian Weyl manifold, the foliation $\D$ is 
Riemannian, totally geodesic. Its leaves, if compact, are complex 
Hopf surfaces (non-primary, in general) admitting an integrable 
hypercomplex structure.
\end{co}
\begin{proof}
Only the second statement has to be proved. It is clear that 
the leaves inherit a hyperhermitian Weyl, non hyperk\"ahler (because $\theta\neq 0$) 
structure. The  
compact hyperhermitian surfaces are classified in \cite{Bo} and the only class 
having the stated property is that of Hopf surfaces.
\end{proof}

As above, here {\it integrable hypercomplex structure} is intended in the sense of
$G$-structures, i.e. of the existence of a local quaternionic coordinate such that the
differential of the change of coordinate belongs to ${\HH}^*$. For further
use we  recall the following:
\begin{te} {\em (cf. \cite{Katt})}\label{lista}
A complex Hopf surface $S$ admits an
integrable hypercomplex structure if and only if $S = ({\HH}
-\{0\})/\Gamma$  where the discrete group $\Gamma$ is conjugate in
$\mathrm{GL}(2, {\CC})$ to any of the following subgroups $G \subset {\HH}^*
\subset  \mathrm{GL}(2, {\CC})$:
 \smallskip
\par $(i)$
$\; G={\ZZ}_m\times \Gamma_c$ with ${\ZZ}_m$ and $\Gamma_c$ both
cyclic generated by left multiplication by $a_m=e^{2\pi i/m}$, $m \geq
1$, and $c\in{\CC}^*$.
\smallskip
\par $(ii)$ $\; G=L\times \Gamma_c$, where $c\in{\RR}^*$ and
$L$ is one of the following: $D_{4m}$,
the dihedral group, $m\geq2$, generated
by the quaternion $j$ and $\rho_m = e^{\pi i/m}$;
$\; T_{24}$, the tetrahedral
group generated by $\zeta^2$ and ${1/ \sqrt 2}(\zeta^3 +\zeta^3j), ~\zeta =
e^{\pi i/4}$;~ $O_{48}$, the octahedral group generated by $\zeta$ and
${1/ \sqrt 2}(\zeta^3+\zeta^3j)$; ~$I_{120}$, the icosahedral group generated
by $\epsilon^3,~j, ~{1/ \sqrt 5}[\epsilon^4 - \epsilon +  (\epsilon^2 -
\epsilon^3)j], ~\epsilon = e^{2\pi i/ 5}.$
\smallskip
\par $(iii)$ ~$G$ generated by ${\ZZ}_m$ and $cj, m\geq 3,~ c\in{\RR}^*$.
\smallskip
\par $(iv)$ ~$G$ generated by $D_{4m}$ and $c\rho_{2n}, ~c\in {\RR}^*$ or by
 $T_{24}$ and $c\zeta, ~ c\in {\RR}^*$.
\end{te}
Contrary to $\D$, the distribution $\D^\perp$ is not integrable. In fact, it plays 
the part of the contact distribution from contact geometry:
\begin{pr}\cite{OP1}\label{non}
On any compact quaternion Hermitian Weyl or hyperhermitian Weyl manifold, 
the distribution $\D^\perp$ is not integrable. Moreover, its integral 
manifolds are totally real and have maximal dimension $n-1$.
\end{pr}
\begin{proof}
Note that a submanifold $N$ is an integral manifold of $\D^\perp$ if and 
only if $\theta$ and $\theta_\al$ vanish on $N$. In this case, also 
$d\theta_\al$ vanishes on $N$. Then \eqref{f} implies that $I_\al X$ is 
normal to $N$ for any $X$ tangent to $N$, \emph{i.e.} $N$ is totally real.
The statement about the dimension of $N$ is now obvious.
\end{proof}
\begin{ex}
Examples of hyperhermitian Weyl manifolds having as leaves of $\D$ any of 
the groups of the surfaces in the above list can be obtained as follows: start 
with the standard hypercomplex Hopf manifold 
$S^1\times S^{4n-1}=\HH^n-\{0\}/\Gamma_2$ (see example \ref{ex}). Consider now the 
diagonal action of any $G$ in Kato's list on $\HH^n$. The action is induced on 
the fibers of the projection $S^\times S^{4n-1}\rightarrow \HH P^{n-1}$, 
hence on the primary standard Hopf surface $S^1\times S^3$ obtaining the desired 
examples. 
\end{ex}

\subsection{Structure theorems}
        In this  section we use the properties of the foliations described 
above to clarify the structure of compact quaternion Hermitian Weyl and 
hyperhermitian Weyl manifolds whose foliations have compact leaves, in 
relation with the other geometries involved: K\"ahler, quaternion K\"ahler, 
locally and globally $3$-Sasakian.

\subsubsection{The link with (locally) $3$-Sasakian geometry}

\begin{te}\cite{OP2}\label{globb}
The class of compact quaternion Hermitian Weyl manifolds $M$ which are not 
quaternion K\"ahler and whose Lee field is quasi-regular, (i.e. 
each point of $M$ has a cubic neighbourhood in which the orbit of $T$ 
enters a finite number of times), 
coincides with the class of flat principal $S^1$-bundles over compact locally 
$3$-Sasakian orbifolds $N=M/T$.
\end{te}
\begin{proof}
Let first $(M,g,H)$ be a compact quaternion Hermitian Weyl manifold as in 
the statement. The orbits of $T$ are closed, hence after rescaling, one 
may suppose they are circles $S^1$ acting on $M$ by isometries because 
$T$ is Killing. The quotient space $N=M/T$ 
is an orbifold (a manifold if $T$ is regular) and, with respect to the 
induced metric $h$, the natural projection $\pi$ becomes a Riemannian 
submersion. Hence, for any leaf $N'$ of $\F$, $\pi_{|N'}:N'\rightarrow N$ is  
a Riemannian covering map. As, according to Proposition \ref{ff}, 
the leaves of $\F$ have a locally $3$-Sasakian structure, $(N,h)$ is locally 
$3$-Sasakian. 

        Conversely, consider a flat principal $S^1$-bundle $\pi:M\rightarrow N$ 
over a compact 
locally $3$-Sasakian manifold $(N,h)$ with local Killing field $\xi_\al$. Choose a  closed $1$-form 
$\theta$ on $M$ defining the flat connection of the bundle $\pi$ and define 
the metric $g:=\pi^*h+\theta\otimes\theta$. Also, define an almost 
quaternionic bundle $H$ on $M$ by defining its local basis as:
\begin{equation}\label{j}
\begin{split}
I_\al=&-\f_\al-\xi_\al^\flat\otimes T, \quad \text{on horizontal fields}\\
I_\al T=&\xi_\al
\end{split}
\end{equation}
where $\f_\al=\nabla^h\xi_\al$ and $T=\theta^\sharp$. 
It is straightforward to check, as in the complex case 
(see \cite{DO}, chapter 6)  
that $(M,g,H)$ is quaternion Hermitian Weyl with Lee form $\theta$.
\end{proof}
\begin{co}\cite{OP1}\label{hom}
The class of compact hyperhermitian Weyl manifolds, not hyperk\"ahler and 
having a quasi-regular (resp. regular) $T$ coincide with the class of flat 
principal $S^1$-bundles over compact $3$-Sasakian orbifolds (resp. manifolds). 
\end{co}

\subsubsection{The link with quaternion K\"ahler geometry}
We now describe the leaf space of the foliation $\D$, when it exists. 
\begin{te}\cite{OP1}, \cite{PPS}\label{ts2}
Let $(M,g,H)$ be a compact quaternion Hermitian Weyl (resp. hyperhermitian Weyl) 
manifold, non 
quaternion K\"ahler (resp. non hyperk\"ahler) whose foliation $\D$ has 
compact leaves. Then the 
leaves space $P=M/\D$ is a compact quaternion K\"ahler orbifold with 
positive scalar curvature, 
the projection is a Riemannian, totally geodesic submersion and a fibre 
bundle map with fibres as described in Proposition \ref{qq} (resp. \ref{qqq}). 
\end{te}
\begin{proof}
In the local case of quaternion Hermitian Weyl $M$, we have to 
explain how to project the structure of $M$ over $P$. The 
key point is that locally, $H$ has admissible basis formed by 
$\nabla$-parallel (hence integrable) complex structures. Then formulae 
\eqref{a}, \eqref{c} show that $H$ is projectable. The foliation being 
Riemannian, $g$ is also projectable. The compatibility of the projected 
quaternion bundle with the projected metric is clear. To show that the projected 
structure is quaternion K\"ahler, let $\omega_P$ be the $4$-form of 
the projected structure. As the projection is a totally geodesic 
Riemannian submersion, $\omega_P$ coincides with the restriction 
of $\omega$ to basic vector fields on $M$. Hence, it is 
enough to show that $\nabla\omega=0$ on basic vector fields. 
But $\nabla\omega=\sum_\al\nabla\omega_\al\wedge\omega_\al+
\omega_\al\wedge\nabla\omega_\al$ 
and the result follows from equation \eqref{b}. The scalar curvature of $(P,g)$ is 
easily computed using O'Neill formulae.

        The global case of a hyperhermitian Weyl $M$ now follows. 

\end{proof}
\begin{re}
The above fibration can never be trivial, according to Proposition \ref{non}.
\end{re}
 
 Let now $M$ be hyperhermitian Weyl, $\T$ be the foliation generated by the 
vector field $T$ and $\V$ the 
$2$-dimensional foliation generated by $T$ and $JT$, where $J$ is a fixed  
compatible 
global complex structure  belonging to $H$. 
Theorem  \ref{ts2}, together with the structure of  
$3$-Sasakian manifolds described in section \ref{3ss}, furnish the following 
structure theorem:
\begin{te}\cite{OP1}, \cite{OP2}\label{diag}  
Let $(M,g,H)$ be a compact  hyperhermitian Weyl 
manifold, non 
 hyperk\"ahler, such that  
the foliations $\D$, $\V$, $\T$ and  $\KK$ have compact leaves. There exists 
the following commutative diagram of fibre bundles and Riemannian submersions 
in the category of orbifolds:

\medskip
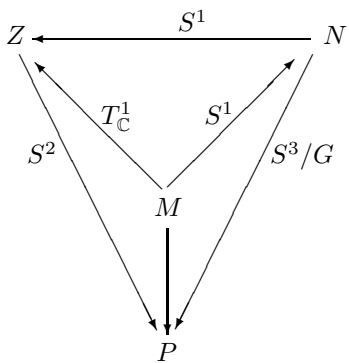
\begin{center}
\begin{picture}(160,140)(10,-40)

\put(67,-40){$P$}
\put(10,80){$Z$}
\put(15,76){\vector(1,-2){52}}
\put(18,35){$S^2$}
\put(130,80){$N$}
\put(125,82){\vector(-1,0){105}}
\put(75,85){$S^1$}
\put(126,76){\vector(-1,-2){52}}
\put(110,35){$S^3/G$}
\put(66,15){$M$}
\put(71,10){\vector(0,-1){40}}
\put(69,25){\vector(-1,1){48}}
\put(71,26){\vector(1,1){48}}
\put(85,50){$S^1$}
\put(46,50){$T^1_\CC$}
\end{picture}
\end{center}

\medskip
Here $N$ is  globally $3$-Sasakian. 
The fibres of $M \rightarrow P$ are Kato's integrable hypercomplex Hopf
surfaces $(S^1 \times S^3)/G$, non necessarily primary and non
necessarily all homeomorphic if $M$ is hyperhermitian Weyl. The $S^1$-bundle 
$P\rightarrow Z$ is a Boothby-Wang fibration.
\end{te}

Note that all arrows appearing in the diagram are canonical,
except for $M \rightarrow Z$, which depends on the choice of
the compatible global complex structure on $M$. However, different choices of 
this complex structure produce analytically equivalent complex manifolds 
$Z$. 
\begin{re}
 The 
diagram \ref{diag} holds also if $dim(M) = 8$. In this case $P$ is
still Einstein by the
above discussion. The integrability of the complex structure on its twistor
space implies it is also self-dual (cf. \cite{Be}).
Then just recall that a
$4$-dimensional $N$ is usually defined to be quaternionic K\"ahler  if it is Einstein
and self-dual.
\end{re}
\begin{re}
For the hyperhermitian Weyl manifold $M=S^1\times S^{4n-1}$, diagram \ref{diag} 
becomes the well-known:

\medskip
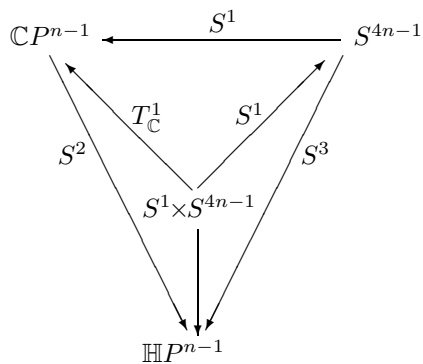
\begin{center}
\begin{picture}(160,140)(0,-40)

\put(50,-40){$\HH P^{n-1}$}
\put(0,80){$\CC P^{n-1}$}
\put(15,76){\vector(1,-2){52}}
\put(18,35){$S^2$}
\put(130,80){$S^{4n-1}$}
\put(125,82){\vector(-1,0){90}}
\put(75,85){$S^1$}
\put(126,76){\vector(-1,-2){52}}
\put(110,35){$S^3$}
\put(50,15){$S^1\!\!\times\!\! S^{4n-1}$}
\put(71,10){\vector(0,-1){40}}
\put(69,25){\vector(-1,1){48}}
\put(71,26){\vector(1,1){48}}
\put(85,50){$S^1$}
\put(46,50){$T^1_\CC$}
\end{picture}
\end{center}
which was the model for the general one. Also, examples of quaternion Hermitian 
Weyl manifolds will be obtained by considering appropriate quotients of 
the manifolds in the vertices of this diagram.
\end{re}

\begin{re}
It is proved in \cite{BGM0} that in every
dimension $4k-5, k\geq 3$ there are infinitely many distinct homotopy
types of complete inhomogeneous 3-Sasakian manifolds. Thus, by simply making
the product with $S^1$, we obtain infinitely many non-homotopically equivalent
examples of compact hyperhermitian Weyl manifolds.
\end{re}

\subsubsection{Some topological consequences of diagram \ref{d}}
A first consequence of the diagram \ref{diag} 
concerns cohomology. Note first that the property $\nabla\theta = 0$ implies
the
vanishing of the Euler characteristic of $M$. 
Then, applying twice the
Gysin sequence in the upper triangle one finds the
relations between the Betti
numbers of $M$ and $Z$ :
\begin{gather*}
b_i(M) = b_i(Z) + b_{i-1}(Z) - b_{i-2}(Z) - b_{i-3}(Z)\hspace{.1 in}
(0 \leq i \leq 2n-1),\\
b_{2n}(M) = 2\left [ b_{2n-1}(Z) - b_{2n-3}(Z)\right ].
\end{gather*}
On the other hand, since $P$ has
positive scalar curvature, both $P$ and its twistor space $Z$ have zero odd
Betti numbers, cf.
\cite{Be}. The Gysin sequence of the fibration $Z \rightarrow P$
then yields:
\begin{equation*}
b_{2p}(Z) = b_{2p}(P) + b_{2p-2}(P)
\end{equation*}
Together with the previous found relations this implies:
\begin{te}\cite{OP1}, \cite{OP2}
Let $M$ be a compact  hyperhermitian Weyl 
manifold satisfying the assumptions of Theorem \ref{diag}. 
Then the following relations hold good:
\begin{gather*}
b_{2p}(M) = b_{2p+1}(M) = b_{2p}(P) - b_{2p-4}(P)\hspace{.1 in}  (0 \leq 2p
\leq 2n-2),\\
b_{2n}(M) =0,\\
\sum_{k=1}^{n-1}k(n-k+1)(n-2k+1)b_{2k}(M)=0.
\end{gather*}
(Poincar\'e duality gives the correspondent  of the first two equalities
for $2n+2 \leq 2p \leq
4n$). In particular $b_1(M) = 1$. Moreover, if $n$ is
even, $M$ cannot carry any quaternion K\"ahler  metric.
\end{te}
The last identity is obtained, by applying S. Salamon's constraints on compact positive quaternion
K\"ahler manifolds to the same diagram (cf \cite{GS}).
\begin{re}
We obtain in particular $b_{2p-4}(P) \leq b_{2p}(P)$ for $0 \leq 2p
\leq 2n-2.$ Since any compact quaternion K\"ahler  $P$ with 
positive scalar curvature can
be realized as the quaternion K\"ahler  base of a compact 
quaternion Hermitian Weyl manifold  $M$, this implies, in the positive scalar
curvature case, the Kraines - Bonan inequalities for
Betti numbers of compact quaternion K\"ahler  manifolds (cf. \cite{Be}).
  
        $b_1(M) = 1$ is a much stronger restriction on the
topology of compact  quaternion Hermitian Weyl manifolds in the larger 
class of compact complex  Vaisman 
(generalized Hopf) manifolds. For the latter,  the only restriction is $b_1$ odd
and the induced Hopf bundles over compact Riemann surfaces of genus $g$
provide examples of Vaisman (generalized Hopf) manifolds with $b_1 = 2g+1$ for any
$g$ , cf. \cite{Va3}.
\end{re}

The properties $b_1 = 1$ and $b_{2n} = 0$ have
the following consequences: 
\begin{co}
Let $(M, I_1, I_2, I_3)$ be a compact hypercomplex manifold that admits
a locally and non globally conformal hyperK\"ahler metric. Then none of the
compatible complex structures $J = a_1I_1 + a_2I_2 + a_3I_3$, $a_1^2 + a_2^2 +
a_3^2 = 1$, can support a K\"ahler  metric. In particular, $(M, I_1, I_2, I_3)$ does
not admit any hyperK\"ahler metric. 

        Let $M$ be a $4n$-dimensional $\mathcal C^\infty$ manifold that 
admits a locally and non globally conformal hyperK\"ahler structure 
$(I_1, I_2, I_3,g)$. Then, for $n$ even, $M$ cannot admit any quaternion K\"ahler  structure and,
for $n$ odd, any quaternion K\"ahler  structure of positive scalar curvature.
\end{co}

\subsubsection{Homogeneous compact hyperhermitian Weyl manifolds}

In the complex case, a complete classification of compact 
homogeneous Vaisman manifolds is still lacking. By contrast, for 
compact homogeneous hyperhermitian Weyl manifolds 
a precise classification may be obtained.
\begin{de}
A hyperhermitian Weyl  manifold $(M,[g],H,D)$ is homogeneous if
there exists a Lie group which acts transitively and effectively on the left
on $M$ by hypercomplex isometries. 
\end{de}
The homogeneity implies the regularity of the canonical foliations:
\begin{te}\cite{OP1}
On a compact 
homogeneous  hyperhermitian Weyl manifold the foliations
$\mathcal D$, $\mathcal V$ and $\mathcal B$  are regular
and in the diagram \ref{diag}, $N$, $Z$, $P$ are
homogeneous manifolds, compatible with the respective structures.
\end{te}
\begin{proof} Fix $J \in H$ be a compatible complex structure on $M$. Then
$(M,g,J)$  is a 
homogeneous Vaisman manifold and by Theorem 3.2 in \cite{Va4} 
we have the
regularity of both the foliations $\mathcal V_J$
and $\mathcal B$. Therefore $M$ projects on homogeneous manifolds $Z_J$ and 
$N$.
In particular the
projections of $I_{\alpha}B$ on $N$ are regular Killing vector fields. Then
Lemma 11.2 in \cite{Tan} assures that the 3-dimensional
foliation spanned by the projections of $I_1B, I_2B, I_3B$ is
regular. This, in turn, implies that $P$ is a homogeneous
manifold, thus $\mathcal D$ is
regular on $M$.
\end{proof}

On the other hand, a compact homogeneous $3$-Sasakian manifolds have 
been classified in \cite{BGM0}. We use this classification together 
with Corollary \ref{hom} to derive:
\begin{pr}\cite{OP1} The class of compact
 homogeneous hyperhermitian Weyl manifolds coincides with that
of flat principal
$S^1$-bundles over one of the 3-Sasakian homogeneous manifolds:
$S^{4n-1}$,
${\RR}P^{4n-1}$ the flag manifolds $\mathrm{SU}(m)/\mathrm{S}(\mathrm{U}(m-2)\times \mathrm{U}(1)),
m\geq 3$, $\mathrm{SO}(k)/$$(\mathrm{SO}(k-4)\times \mathrm{Sp}(1)), k\geq 7$, the exceptional spaces
$G_2/\mathrm{Sp}(1)$, $F_4/\mathrm{Sp}(3)$,
$E_6/\mathrm{SU}(6)$, $E_7/\mathrm{Spin}(12)$,
$E_8/E_7$.
\end{pr}

The flat principal $S^1$-bundles over $P$ are characterized by having
zero or torsion Chern class $c_1 \in H^2(P;\ZZ)$ and classified by it. The
integral cohomology group
$H^2$ of the 3-Sasakian homogeneous manifolds can be computed by
looking at the long homotopy exact sequence
$$...\rightarrow \pi_2(K) \rightarrow \pi_2(G) \rightarrow \pi_2(G/K)
\rightarrow \pi_1(K)\rightarrow \pi_1(G) \rightarrow ... $$
for the 3-Sasakian homogeneous manifolds $G/K$ listed
above. Since $\pi_2(G) = 0$ for any compact Lie group
$G$, one obtains the following isomorphisms (cf. \cite{BGM2}):
\begin{equation*}
H^2~(\frac {\mathrm{SU}(m)}{\mathrm{S}(\mathrm{U}(m-2)\times \mathrm{U}(1))}) \cong \ZZ,
\qquad H^2({\RR}P^{4n-1}) \cong \ZZ_2
\end{equation*}
and $H^2(G/K) = 0$ for
all the other 3-Sasakian homogeneous manifolds. Hence:
\begin{co}\cite{OP1}
Let $M$ be a compact 
 homogeneous hyperhermitian Weyl manifold. Then $M$ is one of the
following:
\par $(i)$ A product $(G/K) \times S^1$, where $G/K$ can be
any of the 3-Sasakian homogeneous manifolds in the list: \\$S^{4n-1}$, 
$\RR P^{4n-1}$,
$\mathrm{SU}(m)/\mathrm{S}(\mathrm{U}(m-2)\times \mathrm{U}(1)), m\geq 3$, $\mathrm{SO}(k)/(\mathrm{SO}(k-4)\times \mathrm{Sp}(1)), k\geq 7$,  $G_2/\mathrm{Sp}(1)$,
$F_4/\mathrm{Sp}(3)$, $E_6/\mathrm{SU}(6)$, $E_7/\mathrm{Spin}(12)$, $E_8/E_7$.
\par $(ii)$ The M\"obius band,
i.e. the unique non  trivial principal $S^1$-bundle over ${\RR}P^{4n-1}$.
\end{co}

        For example in dimension $8$ one obtains only the following spaces: 
$S^7\times S^1$, ${\RR}P^7 \times S^1$, 
$\{\mathrm{SU}(3)/\mathrm{S}(\mathrm{U}(1)\times \mathrm{U}(1))\}\times S^1$ and the M\"obius band over 
${\RR}P^7$. The first exceptional
example appears in  dimension $12$:  the trivial bundle $\{G^2/\mathrm{Sp}(1)\}
\times S^1$ whose $3$-Sasakian base is diffeomorphic to the Stiefel
manifold $V_2({\RR}^7)$ of the orthonormal 2-frames in ${\RR}^7$.
\subsubsection{A hyperhermitian Weyl finite covering of a quaternion 
Hermitian Weyl manifold}

        In general, quaternion K\"ahler manifolds are not finitely covered by 
non simply connected hyperk\"ahler ones. But in the locally conformal K\"ahler 
case we have:
\begin{te}\cite{OP2}
Let $M$ be a compact quaternion Hermitian Weyl manifold which is not quaternion
K\"ahler. If the leaves of ${\mathcal T}$ are compact, then $M$ admits a
finite covering space carrying a structure of a hyperhermitian Weyl manifold.
\end{te}
\begin{proof}
Let first $T$ be a regular vector field. Accordingly, $N=M/T$ is 
compact locally $3$-Sasakian \emph{manifold}, Einstein with positive 
scalar curvature. From Myers theorem, its Riemannian universal cover 
$\tilde N$ is compact and $\pi_1(N)$ is finite. Hence, the pull-back 
$\tilde \KK\rightarrow \tilde N$ (see Corollary \ref{pul}) is trivial 
and $\tilde N$ is globally $3$-Sasakian. Let now $\tilde M\rightarrow \tilde N$ 
be the pull-back of the $S^1$-bundle $M\rightarrow N$: being a 
flat principal circle bundle over a $3$-Sasakian manifold, Corollary 
\ref{hom} provides a hyperhermitian Weyl structure on $\tilde M$. By 
construction, this one projects on the quaternion Hermitian Weyl structure 
of $M$. 

        In the weaker assumption that $\T$ has only compact leaves (it 
is a quasi-regular foliation), the leaves space $N$ is a compact 
orbifold with same Riemannian properties as above. Its universal 
orbifold covering $\tilde N^{orb}$ is a complete Riemannian orbifold with 
positive Ricci curvature. According to Corollary 21 in 
\cite{Borz}, the diameter of 
$\tilde N^{orb}$ is finite. Hence $\tilde N^{orb}$ is compact and 
$\pi_1^{orb}(N)$ is finite. Now the pull-back of $\KK\rightarrow N$ to 
$\tilde N^{orb}$ is again trivial and, as in the manifold case, one shows that 
$\tilde N^{orb}$ is a globally $3$-Sasakian orbifold. The proof then continues 
as above. Note that the total space $\tilde M$ is again a \emph{manifold}.
\end{proof}

\subsection{Examples}
Using the structure theorems, we can now describe a large class of examples 
of quaternion Hermitian Weyl manifolds.

        Recall first that a real $4$-dimensional Hopf manifold is an integrable 
quaternion Hopf manifold, \emph{i.e.} a quotient $(\HH-\{0\})/G=
(\RR^4-\{0\})/G$, 
where $G$ is a discrete subgroup of $\mathrm{CO}(4)\sim \mathrm{GL}(1,\HH)\cdot \mathrm{Sp}(1)$. The 
metric $(h\ov{h})^{-1}dh\otimes d\ov{h}$, globally conformal with the 
flat one on $\HH$, is invariant w.r.t. the action of $G$. This proves: 
\begin{pr}\cite{OP1}
Any real $4$-dimensional Hopf manifold is a compact quaternion Hermitian Weyl 
manifold.
\end{pr}
We generalize this construction to higher dimensions by considering the quaternion 
Hopf manifold $M=(\HH^n-\{0\})/G$, with $G$ of the form \eqref{54}, 
acting diagonally on the 
quaternionic coordiantes $(h^1,...,h^n)$. The metric on $M$ will now  
be the projection of $(\sum_ih^i\ov{h}^i)^{-1}\sum_idh^i\otimes d\ov{h}^i$ and 
is denoted with $g$. 
Moreover, we shall assume the resulting $4$-dimensional foliation  $\D$ to 
have compact leaves. We may state:
\begin{pr}\cite{OP1}
The quaternion Hopf manifold $M=(\HH^n-\{0\})/G$  endowed with the metric $g$ 
is a compact quaternion Hermitian Weyl manifold. The leaves of the foliation 
$\D$ are integrable quaternion Hopf $4$-manifolds. The leaf space  
$P=M/\D$ is a quaternion K\"ahler orbifold quotient of $\HH P^{n-1}$ whose set 
of singular points is, generally, $\RR P^{n-1}\subset \HH P^{n-1}$. Moreover:

         If $G$ is one of the groups in Kato's list 
(see Theorem \ref{lista}), then $M$ is hyperhermitian Weyl, The leaves of $\D$ are 
integrable Hopf surfaces and $P$ is $\HH P^{n-1}$.
\end{pr}
The result follows from the fact that the group $G$, being a discrete subgroup of 
$\mathrm{GL}(n,\HH)\cdot \mathrm{Sp}(1)$,  preserves the quaternionic structure of the universal 
covering of $M$. The structure of the leaves was discussed in Proposition \ref{qq}. 
Note that  $\mathrm{GL}(n,\HH)$ acts on the left and $\mathrm{Sp}(1)$ acts on the right on the 
quaternionic coordinates, hence the induced action of $G$ on 
$\HH P^{n-1}$ fixes the points which can be represented in real coordinates. 
If $G$ belongs to Kato's list, then it is a subgroup of $\mathrm{GL}(n,\HH)$ and 
preserves the hyperhermitian structure of the covering, inducing the 
same structure on the leaves. 

\begin{ex}\cite{PPS}, \cite{OP1}
For $n=2$, let $G$ be the cyclic group generated by $(h^0,h^1)\mapsto 
(2e^{2\pi i/3}h^0, 2e^{4\pi i/3}h^1)$ and $M=(\HH^2-\{0\})/G$. Here the 
leaf space $P=M/\D$ is a $\ZZ_3$ quotient $\HH P^1$. The leaves of $\D$ are 
standard Hopf surfaces $S^1\times S^3$ over the regular points of the orbifold $P$ 
and are non-primary Hopf surfaces  $(S^1\times S^3)/\ZZ_3$ over the two singular 
points of homogeneous coordinates $[1:0]$ and $[0:1]$ of $P$.
\end{ex}

\end{document}